\documentclass[a4paper,11pt,reqno]{elsarticle}
\usepackage{latexsym, amsmath, amsfonts, amssymb, amsthm, amscd, epsfig}
\usepackage{url}

\usepackage[font=small, labelfont={sc}]{caption}
\usepackage{subfig}

\usepackage[a4paper,scale={0.73,0.75},marginratio={1:1},footskip=10mm,headsep=10mm]{geometry}

\usepackage{color}

\usepackage{hyperref}

\numberwithin{equation}{section}

\def \beq {\begin{eqnarray}}
\def \eeq {\end{eqnarray}}
\def \beqn {\begin{eqnarray*}}
\def \eeqn {\end{eqnarray*}}

\newtheorem{theorem}{Theorem}
\newtheorem{itlemma}[theorem]{Lemma}
\newtheorem{itproposition}[theorem]{Proposition}
\newtheorem{itcorollary}[theorem]{Corollary}
\newtheorem{itremark}[theorem]{Remark}
\newtheorem{itdefinition}[theorem]{Definition}
\newtheorem{itexample}[theorem]{Example}
\newtheorem{itclaim}[theorem]{Claim}
\newtheorem{itfact}[theorem]{Fact}

\newenvironment{fact}{\begin{itfact}\rm}{\end{itfact}}
\newenvironment{claim}{\begin{itclaim}\rm}{\end{itclaim}}
\newenvironment{lemma}{\begin{itlemma}}{\end{itlemma}}
\newenvironment{remark}{\begin{itremark}\rm}{\end{itremark}}
\newenvironment{corollary}{\begin{itcorollary}}{\end{itcorollary}}
\newenvironment{proposition}{\begin{itproposition}}{\end{itproposition}}
\newenvironment{definition}{\begin{itdefinition}\rm}{\end{itdefinition}}
\newenvironment{example}{\begin{itexample}\rm}{\end{itexample}}

\newcommand{\be}[1]{\begin{equation}\label{#1}}
\newcommand{\ee}{\end{equation}}
\newcommand{\bl}[1]{\begin{lemma}\label{#1}}
\newcommand{\br}[1]{\begin{remark}\label{#1}}
\newcommand{\brs}[1]{\begin{remarks}\label{#1}}
\newcommand{\bt}[1]{\begin{theorem}\label{#1}}
\newcommand{\bd}[1]{\begin{definition}\label{#1}}
\newcommand{\bp}[1]{\begin{proposition}\label{#1}}
\newcommand{\bc}[1]{\begin{corollary}\label{#1}}
\newcommand{\bfact}[1]{\begin{fact}\label{#1}.}
\newcommand{\bex}[1]{\begin{example}\label{#1}.}
\newcommand{\ec}{\end{corollary}}
\newcommand{\efact}{\end{fact}}
\newcommand{\eex}{\end{example}}
\newcommand{\el}{\end{lemma}}
\newcommand{\er}{\end{remark}}
\newcommand{\ers}{\end{remarks}}
\newcommand{\et}{\end{theorem}}
\newcommand{\ed}{\end{definition}}
\newcommand{\ep}{\end{proposition}}
\newcommand{\epr}{\end{proof}}
\newcommand{\bpr}{\begin{proof}}
\newcommand{\bcl}[1]{\begin{claim}\label{#1}}
\newcommand{\ecl}{\end{claim}}

\newcommand{\ecs}{\end{corollary}}
\newcommand{\eers}{\end{exercise}}
\newcommand{\eexs}{\end{example}}
\newcommand{\eems}{\end{example}}
\newcommand{\els}{\end{lemma}}
\newcommand{\eles}{\end{lemmaex}}
\newcommand{\ets}{\end{theorem}}
\newcommand{\eds}{\end{definition}}
\newcommand{\eps}{\end{proposition}}

\newcommand{\bi}{\begin{itemize}}
\newcommand{\ei}{\end{itemize}}
\newcommand{\ben}{\begin{enumerate}}
\newcommand{\een}{\end{enumerate}}

\def\vbar{\mathchoice{\vrule height6.3ptdepth-.5ptwidth.8pt\kern-.8pt}
   {\vrule height6.3ptdepth-.5ptwidth.8pt\kern-.8pt}
   {\vrule height4.1ptdepth-.35ptwidth.6pt\kern-.6pt}
   {\vrule height3.1ptdepth-.25ptwidth.5pt\kern-.5pt}}
\def\fudge{\mathchoice{}{}{\mkern.5mu}{\mkern.8mu}}
\def\bbc#1#2{{\rm \mkern#2mu\vbar\mkern-#2mu#1}}
\def\bbb#1{{\rm I\mkern-3.5mu #1}}
\def\bba#1#2{{\rm #1\mkern-#2mu\fudge #1}}
\def\bb#1{{\count4=`#1 \advance\count4by-64 \ifcase\count4\or\bba A{11.5}\or
   \bbb B\or\bbc C{5}\or\bbb D\or\bbb E\or\bbb F \or\bbc G{5}\or\bbb H\or
   \bbb I\or\bbc J{3}\or\bbb K\or\bbb L \or\bbb M\or\bbb N\or\bbc O{5} \or
   \bbb P\or\bbc Q{5}\or\bbb R\or\bbc S{4.2}\or\bba T{10.5}\or\bbc U{5}\or
   \bba V{12}\or\bba W{16.5}\or\bba X{11}\or\bba Y{11.7}\or\bba Z{7.5}\fi}}

\def \R {{\mathbb R}}

\def \P {{\mathbb P}}
\def \E {{\mathbb E}}
\def \u {{\bf u}}

\def \ra {\rightarrow }

\def \s {y}

\def \ind {{\bf 1}}

\newcommand{\ba}[1]{\addtocounter{for}{1} \begin{eqnarray}\label{#1}}
\newcommand{\ea}{\end{eqnarray}}

\def\sqr#1#2{{\vcenter{\vbox{\hrule height .#2pt
                             \hbox{\vrule width .#2pt height#1pt \kern#1pt
                                   \vrule width .#2pt}
                             \hrule height .#2pt}}}}

\def\pmb#1{\setbox0=\hbox{#1}%
   \kern-.025em\copy0\kern-\wd0
   \kern.05em\copy0\kern-\wd0
   \kern-.025em\raise.0433em\box0 }
\def\sqr#1#2{{\vcenter{\vbox{\hrule height.#2pt
     \hbox{\vrule width.#2pt height#1pt \kern#1pt
   \vrule width.#2pt}\hrule height.#2pt}}}}

\def\e{\epsilon}

\def\e{\epsilon}
\def\s{\sigma}
\def\d{\delta}
\def\l{\lambda}

\def\g{\gamma}

\def\a{\alpha}
\def\b{\beta}

\newenvironment{myenumerate}{%
\begin{list}{\labelenumi}
	{%
	\setlength{\itemsep}{0.4em}%
	\setlength{\topsep}{0.5em}%
	\setlength\leftmargin{2.6em}%
	\setlength\labelwidth{2.15em}%
	\setlength{\labelsep}{0.45em}%
	\usecounter{enumi}%
	}%
	}%
{\end{list}
}

\renewenvironment{enumerate}{
\renewcommand{\theenumi}{\arabic{enumi}}%
\renewcommand{\labelenumi}{{\rm(\theenumi)}}%
\begin{myenumerate}}%
{\end{myenumerate}}

{\end{myenumerate}}

{\end{myenumerate}}

\newenvironment{myitemize}{%
\begin{list}{$\bullet$}%
 	{%
	\setlength{\itemsep}{0.4em}%
	\setlength{\topsep}{0.5em}%
	\setlength\leftmargin{2.6em}%
	\setlength\labelwidth{2.15em}%
	\setlength{\labelsep}{0.45em}%
	}%
	}%
{\end{list}}

\renewenvironment{itemize}{
\begin{myitemize}}%
{\end{myitemize}}

\begin{document}

\title{Multi-scaling of moments in stochastic volatility models}

\author[pdp]{P. Dai Pra \corref{cor}}
\ead{daipra@math.unipd.it}
\cortext[cor]{Corresponding author}

\author[pp]{P. Pigato}
\ead{pigato@math.unipd.it}

\address[pdp]{Dipartimento di Matematica Pura ed Applicata,
	Universit\`a degli Studi di Padova,
	via Trieste 63,
	I-35121 Padova,
	Italy}

\address[pp]{Dipartimento di Matematica Pura ed Applicata,
	Universit\`a degli Studi di Padova,
	via Trieste 63,
	I-35121 Padova,
	Italy \\  and \\  Laboratoire d'Analyse et de Mathématiques  Appliqu\'ees, Universit\'e Paris-Est Marne-la-Vall\'ee 
Cit\'e Descartes
5, boulevard Descartes 
Champs sur Marne 
77454 Marne la Vallée Cedex 2 France}


\date{\today}

\begin{abstract} 
We introduce a class of stochastic volatility models $(X_t)_{t \geq 0}$ for which the absolute moments of the increments exhibit anomalous scaling: $\E\left(|X_{t+h} - X_t|^q \right)$ scales as $h^{q/2}$ for $q < q^*$, but as $h^{A(q)}$ with $A(q) < q/2$ for $q > q^*$, for some threshold $q^*$. This multi-scaling phenomenon is observed  in time series of financial assets. If the dynamics of the volatility is given by a mean-reverting equation driven by a Levy subordinator and the characteristic measure of the Levy process has power law tails, then multi-scaling occurs if and only if the mean reversion is superlinear.

\end{abstract}

\maketitle

\section{Introduction}

The last few decades have seen a considering effort in constructing stochastic dynamics which exhibit some of the peculiar features of many observed time series, such as: heavy tailed distribution, long memory and path discontinuities. In particular, applications to mathematical finance have motivated the use of stochastic differential equations driven by general Levy processes. In this paper we consider a different, though related, pattern which is rather systematically observed in time series of financial assets, that we call {\em multi-scaling of moments} (\cite{cf:Vas,cf:Gha,cf:Gal,DiM, DIM2}). Let $(X_t)_{t \geq 0}$ be a continuous-time martingale, having stationary increments; in financial applications this could be identified with the de-trended {\em log-price} of an asset, or the price with respect to the martingale measure used to price derivatives. We say the multi-scaling of moments occurs if the limit
\be{ms}
\limsup_{h \downarrow 0} \frac{\log \E\left(|X_{t+h} - X_t|^q \right)}{\log h} =: A(q)
\ee
is {\em non-linear} of the set $\{q \geq 1 : |A(q)| < +\infty\}$. More intuitively, \eqref{ms} says that $\E\left(|X_{t+h} - X_t|^q \right)$ scales, in the limit as  $h \downarrow 0$, as $h^{A(q)}$, with $A(q)$ non-linear. In the case $X_t$ is a Brownian martingale (i.e. a stochastic integral w.r.t. a Brownian motion), one would expect $A(q) = \frac{q}{2}$, at least for $q$ sufficiently small. In this case, multi-scaling of moments can be identified with deviations from this {\em diffusive} scaling, occurring for $q$ above a given threshold; this type of multi-scaling is indeed observed in the values of many financial indexes and exchange rates.

A class of stochastic processes that exhibit multi-scaling for a rather arbitrary scaling function $A(q)$ are the so-called {\em multifractal models} (\cite{cf:CalFisBook, cf:CalFis, cf:CalFisMan}). In these models, the process $X_t$ is given as the {\em random time change} of a Brownian motion:
\be{rtc}
X_t := W_{I(t)},
\ee
where $(W_t)_{t \geq 0}$ is a standard Brownian motion, and $I(t)$ is a stochastic process, often taken to be independent of $W_{\cdot}$, with continuous and increasing trajectories, sometimes called {\em trading time}. Modeling financial series through a random time-change of Brownian motion is a classical topic, dating back to Clark \cite{Cl}, and reflects the natural idea that
external information influences the speed at which exchanges take place
in a market. In multi fractal models, the trading time $I(t)$ is a process with {\em non absolutely continuous} trajectories. As a consequence, $X_t$ cannot be written as a {\em stochastic volatility model}, i.e. in the form $dX_t = \s_t dB_t$, for a Brownian motion $B_t$. This makes the analysis of multi fractal models hard in many respects, as the standard tools of Ito's Calculus cannot be applied.

In \cite{acdp} a much simpler process has been constructed which exhibits a bi-scaling behavior: \eqref{ms} hold with a function $A(q)$ which is piecewise linear and the slope $A'(q)$ takes two different values, which suffices to fit most of the cases observed. This process is a stochastic volatility model, although of a rather peculiar type. Besides exhibiting multi-scaling, this model accounts for other relevant {\em stylized facts} in time series of financial indexes, such as the {\em autocorrelation} profile $t \mapsto Cov(|X_{h} - X_0|, |X_{t+h} - X_t|)$ as well as heavy tailed distribution of $X_{h} - X_0$.

The aim of this paper is to analyze multi-scaling in a more general class of stochastic volatility models, namely those of the form $dX_t = \s_t dB_t$, with a volatility process $\s_t$ independent of the Brownian motion $B_t$; these processes are exactly those that can be written in the form \eqref{rtc} with a trading time $I(t)$ independent of $W_{t}$, and with {\em absolutely continuous} trajectories. We devote special attentions to models in which $V_t := \s^2_t$ is a {\em stationary} solution of a stochastic differential equation of the form
\be{sde}
dV_t \, = \, - f(V_t) dt + dL_t,
\ee
for a {\em Levy subordinator} $L_t$ whose characteristic measure has {\em power law tails at infinity}, and  a function $f(\cdot)$ such that a stationary solution exists, and it is unique in law. We first show multi-scaling is {\em not} possible if $f(\cdot)$ has {\em linear growth}. Thus, the heavy tails produced by the Levy process are not sufficient to produce multi-scaling. On the other hand, we show that, if $f(\cdot)$ behaves as  $ C x^{\g}$ as $x \ra +\infty$, with $C>0$ and $\g>1$, then the stochastic volatility process whose volatility is a stationary solution  of \eqref{sde}, exhibits multi-scaling. In this class of models multi-scaling comes from the combination of heavy tails of $L_t$ and superlinear mean reversion; technically speaking, as will be seen later, the key point is that the distribution of $V_t$ has {\em lighter} tails than those of $L_t$.

We remark that the class of processes introduced in \cite{acdp} can be seen as limiting cases of those considered here, with $\g>2$ and the characteristic measure of the Levy process $L_t$ concentrated on $+\infty$. 

The paper is organized as follows. In Section 2 we give some basic facts on stochastic volatility models, and provide some necessary conditions for multi-scaling. Section 3 contains more specific results for models whose volatility is given by \eqref{sde}.

\section{Multiscaling in stochastic volatility models}

We consider a stochastic process $(X_t)_{t \geq 0}$ that can be expressed in the form
\be{svm}
dX_t = \s_t dW_t,
\ee
where $(W_t)_{t \geq 0}$ is a standard Brownian motion, and $(\s_t)_{t \geq 0}$ is  a {\em stationary}, $[0,+\infty)$-valued process, {\em independent} of $(X_t)_{t \geq 0}$, that we will call the {\em volatility process}. We assume the following weak continuity assumption on the volatility process.
\medskip

\noindent
{\bf Assumption A}. As $h \downarrow 0$, the limit
\[
\frac{1}{h} \int_0^h (\s_s - \s_0)^2 ds \ \ra \ 0
\]
holds in probability.

\medskip
We begin with a basic result on the scaling function $A(q)$ defined in \eqref{ms}. It states that under a uniform integrability condition on the {\em integrated squared volatility}, the diffusive scaling holds. Thus a necessary condition for multi-scaling is the loss of this uniform integrability.

\bp{prof:diff}
Assume that, $p > 1$, 
\be{unifint}
\limsup_{h \downarrow 0} \E \left[\left( \frac{1}{h} \int_0^h \s_s^2 ds \right)^{p/2} \right] < +\infty.
 \ee
Then, under Assumption A,  $A(q) = \frac{q}{2}$ for every $q<p$. 
\ep

\bpr
Note that
\[
\frac{X_h - X_0}{\sqrt{h}} = \frac{1}{\sqrt{h}} \int_0^h \s_s dW_s = \int_0^1 \s_{uh} dB^h_u,
\]
where $B^h_u := \frac{1}{\sqrt{h}} W_{hu}$ is also a standard Brownian motion. Thus, $\frac{X_h - X_0}{\sqrt{h}}$ has the same law of $\int_0^1 \s_{uh} dB_u$, where $B$ is {\em any} Brownian motion independent of the volatility process $(\s_t)_{t \geq 0}$. It follows from Assumption A and the isometry property of the stochastic integral, that 
\be{ui1}
\int_0^1 \s_{uh} dB_u \ \ra \ \s_0 B_1
\ee
in $L^2$ and therefore in probability, as $h \downarrow 0$. By \eqref{unifint} and Burkholder-Davis-Gundy inequality (see \cite{pro}),
\[
\E\left[ \left| \int_0^1 \s_{uh} dB_u \right|^p \right] \leq C_p \E \left[ \left(\int_0^1\s_{uh}^2 du \right)^{p/2} \right] =  \E \left[\left( \frac{1}{h} \int_0^h \s_s^2 ds \right)^{p/2} \right],
\]
so the family of random variables $\left\{ \int_0^1 \s_{uh} dB_u: \, h > 0 \right\}$ is bounded in $L^p$. This implies that  the convergence in \eqref{ui1} is also in $L^q$, for every $q<p$. Thus
\[
\E\left[ \left| \frac{X_h - X_0}{\sqrt{h}} \right|^q \right] = \E\left[ \left| \int_0^1 \s_{uh} dB_u \right|^q \right] \ \ra \ \E\left(\s_0^q \right) \E\left[|B_1|^q \right]
\]
as $h \downarrow 0$ (in particular $\E\left(\s_0^q \right) < +\infty$). Taking the logarithms in the limit above, one obtains $A(q) = \frac{q}{2}$.

\epr

\br{rem:decr}
Suppose $1 \leq q < p $. Then $\frac{A(p)}{p} \leq \frac{A(q)}{q}$. This follows immediately from the fact that, for every $h>0$, 
\[
 \frac{\log \E\left(|X_{t+h} - X_t|^q \right)}{q} = \log \|X_{t+h} - X_t \|_q
 \]
is increasing in $q$, so that $\frac{\log \E\left(|X_{t+h} - X_t|^q \right)}{q \log h}$ is decreasing in $q$ for all $0<h<1$.
\er

In what follows, for models of the form \eqref{svm}, we assume the following further conditions.

\noindent
{\bf Assumption B}. $\E\left( \s_0^2 \right) < +\infty$.

\bigskip

Under Assumption B, \eqref{unifint} holds true for $p=2$.
By Proposition \ref{prof:diff} and Remark \ref{rem:decr}, we have that $A(q) = \frac{q}{2}$ for $1 \leq q <2$, while $\frac{q}{2} \geq A(q) \geq -\infty$ for $q \geq 2$. This suggests the following formal definition of multi-scaling.

\bd{def:multiscaling}
Under Assumptions A and B, we say that multi-scaling occurs if $\{q : - \infty < A(q) <\frac{q}{2}\}$ has a nonempty interior.
\ed

In what follows, Assumptions A and B will be assumed implicitely. Note now that, by Burkholder-Davis-Gundy inequality,  there are constant $c_p,C_p$ such that for each $h>0$
\be{bdg}
c_p \E\left[\left( \int_0^h \s^2_t dt \right)^{p/2} \right] \leq \E\left[\left| X_h - X_0 \right|^p \right] = \E\left[ \left|  \int_0^h \s_s dW_s \right|^p \right] \leq C_p \E\left[\left( \int_0^h \s^2_t dt \right)^{p/2} \right].
\ee
Thus, the condition
\[
 \E \left[\left( \frac{1}{h} \int_0^h \s_s^2 ds \right)^{q/2} \right] < +\infty
 \]
 {\em for each} $h>0$ is necessary for $A(q) > -\infty$. 
 Note also that, by Jensen's inequality,
 \be{tight}
 \E \left[\left( \frac{1}{h} \int_0^h \s_s^2 ds \right)^{q/2} \right] \leq \frac{1}{h}  \int_0^h \E \left[\s_s^q  \right] ds= \E \left[\s_0^q\right],
 \ee
 for $q \geq 2$. Thus, whenever $\E \left[\s_0^q\right] < +\infty$, the assumption of Proposiiton \ref{prof:diff} holds.
 
  This remarks, together with Proposition \ref{prof:diff}, yields the following statement.
 \bc{cor:diff}
 A necessary condition for multi-scaling in \eqref{svm} is that there exists $p >2$ such that
 \[
 \E \left[\left( \frac{1}{h} \int_0^h \s_s^2 ds \right)^{p/2} \right] < +\infty
 \]
 {\em for each} $h>0$, but 
 \[
\E \left[\s_0^p\right] = +\infty.
 \]
 \ec
 
 From the result above we derive an alternative necessary condition for multi-scaling, which has sometimes the advantage to be more easily checked in specific models. 
 
 \bc{cor:diff2}
A necessary condition for multi-scaling in \eqref{svm} is that, for some $h>0$, there exists $p \geq 2$ such that 
\[
\E \left[\s_0^p\right] <  +\infty \  \  \  \  \E \left[\sup_{0 \leq t \leq h} \s_t^p \right] = +\infty.
\]
 
\ec
\bpr
Assume multi-scaling holds, and define
\[
q^* := \inf\{ q : \E \left[\s_0^q\right] = +\infty \}.
\]
By Corollary \ref{cor:diff}, $q^* < +\infty$ while, by Assumption B, $q^* \geq 2$. Moreover, by Proposition \ref{prof:diff}, $A(q) = q/2$ for $q<q^*$. Thus, by Definition \ref{def:multiscaling}, $A(q)$ has to be finite for some $q>q^*$; in particular, as observed above,
\[
 \E \left[\left( \frac{1}{h} \int_0^h \s_s^2 ds \right)^{q/2} \right] < +\infty
 \]
 for $h>0$. Consider $l,r$ with $q^* < l < r < q$. Setting $M_h := \sup_{0 \leq t \leq h} \s_t$, we have
 \[
 \frac{1}{h} \int_0^h \s_s^l ds \leq M_h^{l-2} \frac{1}{h} \int_0^h \s_s^2 ds.
 \]
 By stationarity of $\s_t$, and by applying H\"{o}der inequality with conjugate exponents $\frac{r}{2}$ and $\frac{r/2}{r/2 -1}$, we obtain
 \[
 \E\left(\s_0^l \right) \leq \left[ \E\left( M_h^{r \frac{l/2-1}{r/2 -1}} \right) \right]^{1-\frac{2}{r}} \left[ \E \left[ \left(  \frac{1}{h} \int_0^h \s_s^2 ds \right)^{r/2} \right] \right]^{2/r}.
 \]
 Since $l>q^*$, it follows that $ \E\left(\s_0^l \right) = +\infty$. Moreover, being $r<q$, 
 \[
 \E \left[ \left(  \frac{1}{h} \int_0^h \s_s^2 ds \right)^{r/2} \right] < +\infty.
 \]
 Thus, necessarily, 
 \[
  \E\left( M_h^{r \frac{l/2-1}{r/2 -1}} \right) = +\infty.
  \]
 It is easily checked that, choosing $l$ and $q^*$ sufficiently close, one gets
 \[
 \tilde{r} := r \frac{l/2-1}{r/2 -1} < q^*,
 \]
 which implies
 \[
 \E\left( \s_0^{\tilde{r}} \right) < +\infty.
 \]
Setting $p:= \max(\tilde{r}, 2)$, the proof is completed.

\epr

We conclude this section showing a further property of the scaling function $A(q)$

\br{rem:increasing}
Assume that, for each $h>0$, the integrated volatility has moments of all orders, i.e.
\be{finmom}
\E\left[\left( \int_0^h \s^2_t dt \right)^q \right] < +\infty \ \ \ \mbox{ for every $q \geq 1$.}
\ee
 The following argument shows that, under this assumption, $A(q)$ is increasing in $q$. We will see later an example in which the integrated volatility has heavy tails, so it violates \eqref{finmom}, and $A(\cdot)$ is decreasing in an interval.  We begin by observing that, by \eqref{bdg},
\be{increasing1}
A(q) = \limsup_{h \downarrow 0} \frac{ \log  \E\left[\left( \int_0^h \s^2_t dt \right)^{p/2} \right]}{\log h}.
\ee
From this it easily follows that
\be{increasing2}
\liminf_{h \downarrow 0} \frac{\E\left[\left( \int_0^h \s^2_t dt \right)^{p/2} \right]}{h^{\l}} = 0 \ \Longrightarrow \ \l \leq A(q),
\ee
and
\be{increasing3}
\l < A(q) \ \Longrightarrow \ \liminf_{h \downarrow 0} \frac{\E\left[\left( \int_0^h \s^2_t dt \right)^{p/2} \right]}{h^{\l}} = 0.
\ee
Consider $p>q \geq 1$. Moreover, let $\e>0$, and take $l<q$ such that $[A(q) - \e] \frac{q}{l} < A(q)$. Set
\[
a_h := \left( \int_0^h \s^2_t dt \right)^{1/2}.
\]
We now use Young's inequality $\a \b \leq \frac{\a^r}{r} + \frac{\b^{r'}}{r'}$, valid for $\a,\b \geq 0$, $r,r'>0$, $\frac{1}{r} + \frac{1}{r'} = 1$. Choosing $\a = \frac{a_h^l}{h^{A(q) - \e}}$, $\b = a_h^{p-l}$, $r = \frac{q}{l}$, we get
\[
\frac{a_h^p}{h^{A(q) - \e}} \leq \frac{l}{q} \frac{a_h^q}{h^{(A(q) - \e)\frac{q}{l}}} + \frac{q-l}{q} a_h^{q\frac{p-l}{q-l}}.
\]
Taking expectations:
\be{increasing4}
\frac{\E\left[\left( \int_0^h \s^2_t dt \right)^{p/2} \right]}{h^{A(q) - \e}} \leq \frac{l}{q} \frac{\E\left[\left( \int_0^h \s^2_t dt \right)^{q/2} \right]}{h^{(A(q) - \e)\frac{q}{l}}} + \frac{q-l}{q} \E\left[\left( \int_0^h \s^2_t dt \right)^{q\frac{p-l}{2(q-l)}} \right].
\ee
Since $[A(q) - \e] \frac{q}{l} < A(q)$, by \eqref{increasing3}
\be{increasing5}
\liminf_{h \downarrow 0}  \frac{\E\left[\left( \int_0^h \s^2_t dt \right)^{q/2} \right]}{h^{(A(q) - \e)\frac{q}{l}}}  = 0.
\ee
Moreover, 
\be{increasing6}
\lim_{h \downarrow 0} \E\left[\left( \int_0^h \s^2_t dt \right)^{q\frac{p-l}{2(q-l)}} \right] = 0
\ee
by \eqref{finmom} and dominated convergence. It follows from \eqref{increasing4}, \eqref{increasing5} and \eqref{increasing6}, that
\[
\liminf_{h \downarrow 0} \frac{\E\left[\left( \int_0^h \s^2_t dt \right)^{p/2} \right]}{h^{A(q) - \e}} = 0
\]
which, together with \eqref{increasing2}, yields $A(p) \geq A(q) - \e$. Since $\e$ is arbitrary, the conclusion follows.
\er

\section{Superlinear Ornstein-Uhlenbeck volatility}

In this section we devote our attention to a specific class of stochastic volatility models, namely those of the form
\be{NOU}
\begin{split}
dX_t & = \, \s_t dB_t \\ dV_t & = \, - f(V_t) dt + dL_t \\ V_t & = \, \s^2_t,
\end{split}
\ee
where:
\bi
\item
$(B_t)_{t \geq 0}$ is a standard Brownian motion.
\item
 $(L_t)_{t \geq 0}$ is a Levy process with increasing paths ({\em subordinator})  independent of $(B_t)_{t \geq 0}$. More precisely $(L_t)_{t \geq 0}$ is a real-valued process, with independent increments, $L_0 = 0$ and 
 \[
\E\left[ \exp (- \l L_t) \right] = \exp[-t \Psi(\l)],
\]
with
\[
\Psi(\l) = m \l + \int_{(0,+\infty)} \left(1-e^{-\l x}\right) \nu(dx),
\]
where $m \geq 0$ is the {\em drift} of the process, and $\nu$ is a positive measure on $(0,+\infty)$, called {\em characteristic measure}, satisfying the condition
\[
\int_{(0,+\infty)} (1 \wedge x) \nu(dx) < +\infty.
\]
For generalities on Levy Processes see \cite{ber, CoTa, sato}.
\item
$f(\cdot)$ is a locally Lipschitz, nonnegative function such that $f(0) = 0$ (which guarantees $V_t \geq 0$ if $V_0 \geq 0$). 
\ei
Some conditions on $f(\cdot)$ are needed for \eqref{NOU} to have a stationary solution. We will address this point later. We will always assume that $V_0$ is independent of $(L_t)_{t \geq 0}$. We note now that for many ``natural'' choices of $f$, multi-scaling is not allowed. In particular, multiscaling is not present in Ornstein-Uhlenbeck models (see e.g. \cite{kl1, kl2}).
\bp{prop:linOU}
Suppose $f(\cdot)$ satisfies the {\em linear growth condition}
\[
|f(v)| \leq Av +B
\]
for some $A,B>0$ and all $v>0$. Moreover, assume \eqref{NOU} has a solution for which $(V_t)_{t \geq 0}$ is stationary, nonnegative  and integrable, such that Assumptions A and B hold. Then multi-scaling does not occur.
\ep
\bpr
By Remark \ref{rem:decr}, $A(q) \leq q/2$, so we need to show the converse inequality. Let $(V'_t)_{t \geq }$ be solution of
\be{linOU}
\begin{split}
dV'_t & = -(AV'_t + 2B) dt + dL_t \\ V'_0 & = V_0.
\end{split}
\ee
Note that
\[
d(V_t - V'_t) = -\left[ f(V_t)-AV'_t - 2B \right]dt.
\]
In particular $V_t - V'_t$ is continuously differentiable, and $V_0 - V'_0 = 0$. It follows that $V_t - V'_t \geq 0$ for every $t \geq 0$: indeed the path of $V_t - V'_t$ cannot downcross the value zero, since whenever $\overline{t}$ is such that $V_{\overline{t}} = V'_{\overline{t}} = v$, then
\[
\frac{d}{dt} (V_{\overline{t}} - V'_{\overline{t}} ) = - f(v) + Av + 2B \geq B >0.
\]
Thus for every $t \geq 0$
\[
V_t \geq V'_t = V_0 e^{-At} + \frac{2B}{A}\left(e^{-At} -1 \right) + \int_0^t e^{-A(t-s)} dL_s \geq  V_0 e^{-At} + e^{-tA/2} L_{t/2} - \frac{2B}{A}.
\]
On the other hand
\[
V_t = V_0 - \int_0^t f(V_s) ds + L_t \leq V_0 + L_t,
\]
which yields
\[
\sup_{t \in [0,h]} V_t \leq V_0 + L_h.
\]
Putting all together
\be{linOU1}
V_0 e^{-Ah} + e^{-Ah/2} L_{h/2} - \frac{2B}{A} \leq V_h \leq \sup_{t \in [0,h]} V_t \leq V_0 + L_h.
\ee
Since
\[
V_0 e^{-Ah} + e^{-Ah/2} L_{h/2} - \frac{2B}{A} \in L^p \ \iff \ V_0 + L_h \in L^p,
\]
the conclusion now follows from \eqref{linOU1} and Corollary \ref{cor:diff2}.

\epr

Proposition \ref{prop:linOU} shows that, for models of the form \eqref{NOU} to exhibit multi-scaling, one need to consider a drift $f(\cdot)$ with a superlinear growth. 

\bd{def:regvar}
We say that a function $f:(0,+\infty) \ra (0,+\infty)$ is {\em regularly varying at infinity} with exponent $\a \in \R$ if, for every $x>0$,
\[
\lim_{t \ra +\infty} \frac{f(tx)}{f(t)} = x^{\a}.
\]

\ed

In the case $\a=0$ we say that $f$ is {\em slowly varying at infinity}. Note that $f$ is regularly varying at infinity with exponent $\a$ if and only if $f(u) = u^{\a} l(u)$ where $l$ is slowly varying at infinity. In what follows we consider models of the form \eqref{NOU} for which  the following assumptions hold:
\bi
\item[{\bf A1}]
$(B_t)_{t \geq 0}$ is a standard Brownian motion.
\item[{\bf A2}]
$(L_t)_{t \geq 0}$ is a  Levy subordinator with  characteristic measure $\nu$. Moreover $(B_t)_{t \geq 0}$ and $(L_t)_{t \geq 0}$ are independent.
\item[{\bf A3}]
The function $u \mapsto \nu((u,+\infty))$ is regularly varying at infinity with exponent $-\a<0$.
\item[{\bf A4}]
$f:[0,+\infty) \ra [0,+\infty)$ is increasing, locally Lipschitz, $f(0) = 0$, and it is regularly varying at infinity with exponent $\g>1$.
\ei

The following result has been proved in \cite{sam} (see also \cite{KoYa} for related results).
\bt{th:stmeas}
Under assumption A2-A4, the equation $dV_t \, = \, - f(V_t) dt + dL_t$ admits an unique stationary distribution $\mu$. Moreover $\mu((u,+\infty))$ is regularly varying at infinity with exponent $-\a-\g+1$.
\et

In what follows we assume $V_0$ is independent of $(B_t)_{t \geq 0}$ and $(L_t)_{t \geq 0}$, and it has distribution $\mu$. Theorem \ref{th:stmeas} shows that, if $\g>1$, $V_t$ has a distribution with lighter tails 
than those of the Levy process $L_t$.

We are now ready to state the main result of this paper.
\bt{th:main}
Assume A1-A4 are satisfied, and that $\a+\g>2$ (which, in particular, implies Assumption B). Then the following statements hold.
\ben
\item
If $\g \geq 2$ then
\[
A(q) = \left\{ \begin{array}{ll} \frac{q}{2} & \mbox{for }  1 \leq q < 2(\a+\g-1) \\ \frac{\gamma-2}{2(\gamma-1)}q+\frac{\a+\gamma-1}{\gamma-1} & \mbox{for } q> 2(\a+\g-1). \end{array} \right.
\]
\item
If  $1<\g<2$ then
\[
A(q) = \left\{ \begin{array}{ll} \frac{q}{2} & \mbox{for }  1 \leq q < 2(\a+\g-1) \\ \frac{\gamma-2}{2(\gamma-1)}q+\frac{\a+\gamma-1}{\gamma-1} & \mbox{for } 2(\a+\g-1)< q < \frac{2\a}{2-\g} \\ -\infty & \mbox{for } q>\frac{2\a}{2-\g}. \end{array} \right.
\]

Moreover, for $q \neq 2(\a+\g-1), \frac{2\a}{2-\g}$, the scaling exponent $A(q)$  in \eqref{ms} can be defined as a limit rather that a $\limsup$.
\een
\et

We remark that, in the case $1<\g<2$, $A(\cdot)$ is decreasing for $2(\a+\g-1)< q < \frac{2\a}{2-\g}$. This is not in contradiction with Remark \ref{rem:increasing}, since assumption \eqref{finmom}   is not satisfied. 
\br{rem:linear}
A simple consequence of Theorem \ref{th:main}, is that, by a comparison argument, Proposition \ref{prop:linOU} can be extended to any $f$ which is regularly varying at infinity with exponent $1$.
\er

The proof of Theorem \ref{th:main} will be divided into several steps. We begin by dealing with the case $f(v) = C v^{\g}$, with $C>0$, and $L_t$ is a {\em compound Poisson process}.

\bp{prop:main1}
The conclusion of Theorem \ref{th:main} hold if $f(v) = C v^{\g}$, with $C>0$, $L_t$ is a Levy subordinator with zero drift and {\em finite}  characteristic measure $\nu$.
\ep
\bpr
Note that, for $q < 2(\a+\g-1)$, by Theorem \ref{th:stmeas}, we have $E\left[ V_0^{q/2} \right] < +\infty$ so that, by Proposition \ref{prof:diff} and \eqref{tight},  $A(q) = \frac{q}{2}$. Thus it is enough to consider the case $q>2(\a+\g-1)$. In what follows we also write $a_h \sim h^u$ for
\be{logequiv}
\lim_{h \ra 0} \frac{ \log a_h}{\log h} = u.
\ee
We will repeatedly use the simple fact that \eqref{logequiv} follows if we show that for every $\e>0$ there exist $C_{\e} >1$ such that
\[
\frac{1}{C_{\e}} h^{u+\e} < a_h < C_{\e} h^{u-\e}.
\]

In what follows all estimates on $A(q)$ are based on  the fact (see \eqref{increasing1}) that the limit
\[
\lim_{h \downarrow 0} \frac{ \log  \E\left[\left( \int_0^h \s^2_t dt \right)^{p/2} \right]}{\log h}
\]
exists if and only if the limit
\[
\lim_{h \downarrow 0} \frac{\log \E\left(|X_{t+h} - X_t|^q \right)}{\log h} 
\]
exists, and in this case the coincide.

\noindent
\textbf{Part 1: $\gamma>2$}\\
By the assumption of finiteness of $\nu$, $(L_t)$ jumps finitely many times in any compact interval. Denote by $(T_k)_{k \geq 1}$ the (ordered) set of positive jump times, and $T_0=0$. Given $h>0$, we denote by $i(h)$ the random number of jump times in the interval $(0,h]$. 

\noindent
{\em Case $i(h) = 0$}. When $i(h) = 0$, $V_t$ solves, for $t \in [0,h]$, the differential equation $\frac{d}{dt} V_t = - C V_t^{\g}$, whose solution is
\[
V_t=\left( V_0^{1-\gamma}+(\gamma-1)Ct\right)^\frac{1}{1-\gamma}.
\]
Integrating, we get
\be{main1}
\int_0^h V_t dt=\frac{\gamma-2}{\gamma-1}\left[(V_0^{1-\gamma}+(\gamma-1)Ch)^\frac{\gamma-2}{\gamma-1}-(V_0^{1-\gamma})^\frac{\gamma-2}{\gamma-1}\right].
\ee
Note that, setting $\l := \nu([0,+\infty))$,
\[
\E\left[\left(\int_0^h V_t dt\right)^{q/2} \ind_{\{i(h) = 0\}} \right] = \E\left[  \left(\int_0^h V_t dt\right)^{q/2} \Big| i(h)=0 \right] e^{-\l h}.
\]
The factor $e^{-\l h}$ gives no contribution to the behavior of $\E\left[\left(\int_0^h V_t dt\right)^{q/2} \ind_{\{i(h) = 0\}} \right]$ as $h \ra 0$, and it can be neglected. Moreover, by \eqref{main1}, and using the fact that $V_0$ and $\{i(h)=0\}$ are independent,
\begin{multline} \label{main2}
\E\left[\left(\int_0^h V_t dt\right)^{q/2} \Big| i(h)=0  \right]   = \\ 
=  \left(\frac{\g-2}{\g-1}\right)^{q/2} ((\gamma-1)Ch)^{\frac{\gamma-2}{2(\gamma-1)}q}
\E\left[\left[\left(\frac{V_0^{1-\gamma}}{(\gamma-1)Ch}+1\right)^\frac{\gamma-2}{\gamma-1}-\left(\frac{V_0^{1-\gamma}}{(\gamma-1)Ch}\right)^\frac{\gamma-2}{\gamma-1}\right]^{q/2} \right].
\end{multline}
Since, for $0<a<1$ and $z>0$, 
\be{tecineq1}
a(z+1)^{a-1} \leq (z+1)^{a} - z^{a} \leq (z+1)^{a-1},
\ee
for computing the limit $\lim_h \frac{\log \E\left[\left(\int_0^h V_t dt\right)^{q/2} \right]}{\log h}$, the right hand side of \eqref{main2} can be replaced by (using the previous inequality for $a= \frac{\g-2}{\g-1}$; recall that $\g>2$)
\be{main3}
h^{\frac{\gamma-2}{2(\gamma-1)}q}
\E\left[\left[\left(\frac{V_0^{1-\gamma}}{(\gamma-1)Ch}+1\right)^{\frac{\gamma-2}{\gamma-1}-1}\right]^{q/2}\right] 
= h^{\frac{\gamma-2}{2(\gamma-1)}q} \E\left[\left(\frac{V_0^{1-\gamma}}{(\gamma-1)Ch}+1\right)^{-\frac{q}{2(\g-1)}}\right].
\ee
In other words:
\be{main4}
\E\left[\left(\int_0^h V_t dt\right)^{q/2} \ind_{\{i(h) = 0\}} \right] \sim h^{\frac{\gamma-2}{2(\gamma-1)}q} \E\left[\left(\frac{V_0^{1-\gamma}}{(\gamma-1)Ch}+1\right)^{-\frac{q}{2(\g-1)}}\right].
\ee
To estimate the r.h.s. of \eqref{main4}, we observe that for $y>0$ and $0<u<r$, the following inequalities can be easily checked
\be{main5}
\frac{1}{2^r} \ind_{\{y<1\}} \leq (1+y)^{-r} \leq (1+y)^{-u} \leq y^{-u}.
\ee
Setting $r:= \frac{q}{2(\g-1)}$ and $Y := \frac{V_0^{1-\gamma}}{(\gamma-1)Ch}$, using \eqref{main5} we obtain
\be{main6}
\frac{1}{2^r} \P(Y<1) \leq \E\left[\left(\frac{V_0^{1-\gamma}}{(\gamma-1)Ch}+1\right)^{-\frac{q}{2(\g-1)}}\right] \leq E\left(Y^{-u}\right)
\ee
for every $u <  \frac{q}{2(\g-1)}$.
Set $\xi := \frac{\a+\g-1}{\g-1}$. Note that $\xi<r$ for $q>2(\a+\g-1)$. By Theorem \ref{th:stmeas}:
\be{main7}
\P(Y<1) = \P\left(V_0 > \left(\frac{1}{(\g-1)Ch}\right)^{\frac{1}{\g-1}} \right) \sim\left[ \left(\frac{1}{(\g-1)Ch}\right)^{\frac{1}{\g-1}}\right]^{\a+\g-1} \sim h^{\xi},
\ee
Moreover, take $u< \xi$. We have
\be{main8}
E\left(Y^{-u}\right) = \left[(\g-1)Ch\right]^u \E\left[ V_0^{u(\g-1)}\right] \leq A h^u,
\ee
for some $A>0$ that may depend on $u$ but not on $h$,
where we have used the fact that $ \E\left[ V_0^{u(\g-1)}\right] <+\infty$, since $u(\g-1)< \a+\g-1$.
Since $u$ can be taken arbitrarily close to $\xi$, by \eqref{main6}, \eqref{main7} and \eqref{main8} we obtain
\be{main9}
\E\left[\left(\frac{V_0^{1-\gamma}}{(\gamma-1)Ch}+1\right)^{-\frac{q}{2(\g-1)}}\right] \sim h^{\xi},
\ee
which yields
\be{main10}
\E\left[\left(\int_0^h V_t dt\right)^{q/2} \ind_{\{i(h) = 0\}} \right] \sim h^{\frac{\gamma-2}{2(\gamma-1)}q+ \frac{\a+\g-1}{\g-1}}.
\ee
Note that \eqref{main10} has the right order, according to the statement of Theorem \ref{th:main}. Therefore, in order to complete the proof for $\g>2$, it is enough to show that for each $u< \frac{\a+\g-1}{\g-1}$
\be{purpose}
\E\left[\left(\int_0^h V_t dt\right)^{q/2} \ind_{\{i(h) \geq  1 \}} \right] \leq A h^{\frac{\gamma-2}{2(\gamma-1)}q+ u}
\ee
for some $A>0$ that may depend on $u$ but not on $h$.

\medskip
\noindent
{\em Case $i(h) = 1$}. Now
\[
V_t = \left\{ \begin{array}{ll} \left( V_0^{1-\gamma}+(\gamma-1)Ct\right)^\frac{1}{1-\gamma} & \mbox{for } 0 \leq t \leq T_1 \\ \left( V_{T_1}^{1-\gamma}+(\gamma-1)C(t-T_1)\right)^\frac{1}{1-\gamma} & \mbox{for } T_1 \leq t \leq h, \end{array} \right.
\]
which yields
\be{main11}
\begin{split}
\int_0^h V_t dt &=\frac{\gamma-2}{\gamma-1}\left[(V_0^{1-\gamma}+(\gamma-1)C T_1)^\frac{\gamma-2}{\gamma-1}-(V_0^{1-\gamma})^\frac{\gamma-2}{\gamma-1}\right]
\\ 
& + \frac{\gamma-2}{\gamma-1}\left[(V_{T_1}^{1-\gamma}+(\gamma-1)C(h-{T_1}))^\frac{\gamma-2}{\gamma-1}-(V_{T_1}^{1-\gamma})^\frac{\gamma-2}{\gamma-1}\right] \\ & =: P(h) + Q(h),
\end{split}
\ee
and therefore
\be{main12}
\E\left[\left(\int_0^h V_t dt\right)^{q/2} \ind_{\{i(h) = 1\}} \right] \leq 2^{q-1} \left[ \E\left(P^{q/2}(h) \ind_{\{i(h) = 1\}}\right) + \E\left(Q^{q/2}(h) \right)\ind_{\{i(h) = 1\}} \right]
\ee
We now show that $\E\left[\left(\int_0^h V_t dt\right)^{q/2} \ind_{\{i(h) = 1\}} \right] $ can be {\em bounded above} as in \eqref{purpose}:
\be{bound1}
\E\left[\left(\int_0^h V_t dt\right)^{q/2} \ind_{\{i(h) = 1\}} \right] \leq A h^{\frac{\gamma-2}{2(\gamma-1)}q+ u},
\ee
for every $u<\frac{\a+\g-1}{\g-1}$.
By \eqref{main12} it suffices to show that both $\E\left(P^{q/2}(h) \ind_{\{i(h) = 1\}}\right)$ and $\E\left(Q^{q/2}(h) \ind_{\{i(h) = 1\}}\right)$ have an upper bound of the same form. \\
Note first that
\[
P(h) \leq \frac{\gamma-2}{\gamma-1}\left[(V_0^{1-\gamma}+(\gamma-1)C h)^\frac{\gamma-2}{\gamma-1}-(V_0^{1-\gamma})^\frac{\gamma-2}{\gamma-1}\right],
\]
which coincides with \eqref{main1}, whose scaling has already been obtained. Since $\P(i(h)=1) \sim h$, we have that $\E\left(P^{q/2}(h) \ind_{\{i(h) = 1\}}\right)$ scales as the term studied in the case $i(h) = 0$, but with an extra factor $h$, i.e.
\be{main13}
\E\left(P^{q/2}(h) \ind_{\{i(h) = 1\}}\right) \leq A h^{1+\frac{\gamma-2}{2(\gamma-1)}q+ u} \leq A h^{\frac{\gamma-2}{2(\gamma-1)}q+ u}.
\ee
For the term $\E\left(Q^{q/2}(h) \ind_{\{i(h) = 1\}}\right)$ we repeat the steps of the case $i(h) = 0$ (note that all inequalities used there held pointwise) with $V_{T_1}$ in place of $V_0$ and $h-T_1$ in place of $h$ (see \eqref{main4}), obtaining
\be{main14}
\begin{split}
\E\left(Q^{q/2}(h) \ind_{\{i(h)  = 1\}}\right)  & \leq  \E\left[(h-T_1)^{\frac{\gamma-2}{2(\gamma-1)}q}\left(\frac{V_{T_1}^{1-\gamma}}{(\gamma-1)C(h-T_1)}+1\right)^{-\frac{q}{2(\g-1)}} \ind_{\{i(h) = 1\}}\right] \\ & \leq h^{\frac{\gamma-2}{2(\gamma-1)}q}\E\left[\left(\frac{V_{T_1}^{1-\gamma}}{(\gamma-1)C(h-T_1)}+1\right)^{-\frac{q}{2(\g-1)}} \ind_{\{i(h) = 1\}}\right]
\end{split}
\ee
This last term can be bounded from above as follows, for $u< \frac{q}{2(\g-1)}$ and using the trivial bound $V_{T_1} \leq V_0 + L_{h}$
\be{main15}
\begin{split}
 \E\left[\left(\frac{V_{T_1}^{1-\gamma}}{(\gamma-1)C(h-T_1)}+1\right)^{-\frac{q}{2(\g-1)}} \ind_{\{i(h) = 1\}}\right]  & \leq   \E\left[\left(\frac{V_{T_1}^{1-\gamma}}{(\gamma-1)C(h-T_1)}+1\right)^{-u} \ind_{\{i(h) = 1\}}\right] \\
& \leq  \E\left[\left(\frac{V_{T_1}^{1-\gamma}}{(\gamma-1)C(h-T_1)}
\right)^{-\u} \ind_{\{i(h) = 1\}}\right] \\
& \leq  \E\left[\left(\frac{V_{T_1}^{1-\gamma}}{(\gamma-1)Ch)}
\right)^{-\u} \ind_{\{i(h) = 1\}}\right] 
 \\ & \leq A h^{u} \E\left[\left(V_0 + L_{h}\right)^{u(\gamma-1)}
\ind_{\{i(h) = 1\}}\right]
\end{split}
\ee
for a constant $A>0$.
Now observe that $V_0$   is independent of $\ind_{\{i(h) = 1\}}$, that $L_{h}$ has distribution $\nu$ conditioned to $\{i(h) = 1\}$, and that $\P(i(h) = 1) \leq \l h$. It follows that, for a suitable constant $C>0$,
\be{main15bis}
\begin{split}
\E\left[\left(V_0 + L_{h}\right)^{u(\gamma-1)} \ind_{\{i(h) = 1\}}\right] & \leq C P(i(h) = 1) \left[ \E\left(V_0^{u(\gamma-1)}\right) + \E\left(L_{h}^{u(\gamma-1)} | i(h) = 1\right) \right] \\
& = C P(i(h) = 1) \left[ \int v^{u(\gamma-1)} \mu(dv) + \int l^{u(\gamma-1)} \nu(dl) \right].
\end{split}
\ee
Since, by Theorem \ref{th:stmeas}, the tails of $\mu$ are lighter that those of $\nu$, the above integrals are both finite if an only if $\int v^{u(\g-1)} \nu(dv)<+\infty$, which holds true for $u<\a/(\g-1)$ (assumption A3). Thus, for every $u<\a/(\g-1)$,
\be{main16}
\E\left[\left(V_0 + L_{h}\right)^{u(\gamma-1)} \ind_{\{i(h) = 1\}}\right] \leq A h
\ee
for some $A>0$.
 By \eqref{main14}, \eqref{main15}, \eqref{main15bis} and \eqref{main16}, we have that $\E\left(Q^{q/2}(h) \ind_{\{i(h) = 1\}}\right)$ is bounded from above by 
 \[
B h^{\frac{\gamma-2}{2(\gamma-1)}q +u+1}
\]
for every $u<\a/(\g-1)$ and some $B>0$ possibly depending on $u$. Equivalently,
\be{main17}
\E\left(Q^{q/2}(h) \ind_{\{i(h) = 1\}}\right)  \leq B h^{\frac{\gamma-2}{2(\gamma-1)}q+ u}
\ee
for all $u<\frac{\a+\g-1}{\g-1}$.
Therefore, by \eqref{main13} and \eqref{main17},
\eqref{bound1} is established.

\noindent
{\em Case $i(h) \geq 2$}. To prove \eqref{purpose} and thus to complete the whole proof, we are left to show that
\be{purpose2}
\E\left[\left(\int_0^h V_t dt\right)^{q/2} \ind_{\{i(h) \geq  2 \}} \right] \leq A h^{\frac{\gamma-2}{2(\gamma-1)}q+ u}
\ee
for all $u<\frac{\a+\g-1}{\g-1}$ and some $A>0$.

Let $n \geq 2$, and restrict to the event $\{i(h) = n\}$. We have
\[
V_t = \left\{ \begin{array}{cl} \left( V_0^{1-\gamma}+(\gamma-1)Ct\right)^\frac{1}{1-\gamma} & \mbox{for } 0 \leq t \leq T_1 \\ \left( V_{T_1}^{1-\gamma}+(\gamma-1)C(t-T_1)\right)^\frac{1}{1-\gamma} & \mbox{for } T_1 \leq t \leq T_2 \\ \vdots & \\  \left( V_{T_{n-1}}^{1-\gamma}+(\gamma-1)C(t-T_{n-1})\right)^\frac{1}{1-\gamma} & \mbox{for } T_{n-1} \leq t \leq T_n \\  \left( V_{T_n}^{1-\gamma}+(\gamma-1)C(t-T_n)\right)^\frac{1}{1-\gamma} & \mbox{for } T_n \leq t \leq h,
\end{array} \right.
\]
so that \eqref{main11} becomes
\be{main19}
\begin{split}
\int_0^h V_t dt &=\frac{\gamma-2}{\gamma-1}\sum_{k=1}^n\left[(V_{k-1}^{1-\gamma}+(\gamma-1)C (T_k- T_{k-1}))^\frac{\gamma-2}{\gamma-1}-(V_{k-1}^{1-\gamma})^\frac{\gamma-2}{\gamma-1}\right]
\\ 
& + \frac{\gamma-2}{\gamma-1}\left[(V_{T_n}^{1-\gamma}+(\gamma-1)C(h-{T_n}))^\frac{\gamma-2}{\gamma-1}-(V_{T_n}^{1-\gamma})^\frac{\gamma-2}{\gamma-1}\right] \\ & =: \sum_{k=1}^n P_k(h) + P_{n+1}(h).
\end{split}
\ee
Each term $\E\left[P^{q/2}_k(h) \ind_{\{i(h) = n\}}\right]$ can be estimated as in \eqref{main14} and \eqref{main15}, obtaining
\be{main191}
\begin{split}
\E\left[P^{q/2}_k(h) \ind_{\{i(h)  = n\}}\right] & \leq C h^{\frac{\gamma-2}{2(\gamma-1)}q +u }\E\left[\left(V_0 + L_{h}\right)^{u(\gamma-1)} \ind_{\{i(h) = n\}} \right] 
\\
& \leq C' h^{\frac{\gamma-2}{2(\gamma-1)}q +u} \P(i(h) = n) \left[ \E\left(V_0^{u(\gamma-1)}\right) + \E\left(L_{h}^{u(\gamma-1)} | i(h) = n \right) \right]
\end{split} 
\ee
for $u< \frac{q}{2(\g-1)}$ and some constant $C,C'$ that may depend on $u$ but not on $n$ and $h$.
The distribution of $L_h$ {\em given} $\{i(h) = n\}$ is given by the $n$-fold convolution $\nu^{*n}$. In other words, if $X_1,X_2,\ldots,X_n$ are independent random variables with law $\nu$,
\[
\E\left(L_{h}^{u(\gamma-1)} | i(h) = n \right) = \E\left[\left(X_1 + X_2 + \cdots + X_n \right)^{u(\gamma -1)} \right] \leq n^{u(\gamma -1)-1} \E\left[X_1^{u(\gamma -1)} \right].
\]
For $u < \a/(\gamma -1)$, $\E\left[X_1^{u(\gamma -1)} \right] < +\infty$ as well as $\E\left(V_0^{u(\gamma-1)}\right) < +\infty$. Thus
\be{main192}
\E\left[P^{q/2}_k(h) \ind_{\{i(h)  = n\}}\right]  \leq  C h^{\frac{\gamma-2}{2(\gamma-1)}q +u} n^{u(\gamma -1)-1} \P(i(h) = n),
\ee
for some constant $C$ independent of $n$, $h$ and $k$.
By \eqref{main19}, \eqref{main191} and \eqref{main192} we obtain, for $u< \frac{q}{2(\g-1)}$,
\be{main193}
\begin{split}
\E\left[\left(\int_0^h V_t dt\right)^{q/2} \ind_{\{i(h) = n\}} \right] & \leq n^{q/2} \sum_{k=1}^{n+1} \E\left[P^{q/2}_k(h) \ind_{\{i(h)  = n\}}\right] \\ &  \leq C h^{\frac{\gamma-2}{2(\gamma-1)}q +u} n^{u(\gamma -1) + q/2} \P(i(h) = n).
\end{split}
\ee
We can now sum over $n \geq 2$, observing that $ \P(i(h) = n) \leq \frac{\l^n h^n}{n!}$:
\be{main194}
\begin{split}
\sum_{n \geq 2} \E\left[\left(\int_0^h V_t dt\right)^{q/2} \ind_{\{i(h) = n\}} \right] & \leq C h^{\frac{\gamma-2}{2(\gamma-1)}q +u+2} \sum_{n \geq 2} n^{u(\gamma -1) + q/2} \frac{\l^n h^{n-2}}{n!} \\ & \leq C'  h^{\frac{\gamma-2}{2(\gamma-1)}q +u+2}.
\end{split}
\ee
Since $\frac{q}{2(\g-1)} + 2 > \frac{\a+\g-1}{\g-1}$ (recall that $q>2(\a+\g-1)$), 
 we have  that
\[
\E\left[\left(\int_0^h V_t dt\right)^{q/2} \ind_{\{i(h) \geq 2\}} \right]
\]
is negligible with respect to \eqref{main10}.

\noindent
This  completes the proof for the case $\g>2$.

\noindent
\textbf{Part 2: $1<\gamma<2$}\\

\noindent
{\em Case $i(h) = 0$}. Formula \eqref{main1} still hold, but now $\g-2<0$. So \eqref{main2} becomes
\begin{multline} \label{main20}
\E\left[\left(\int_0^h V_t dt\right)^{q/2} \ind_{\{i(h) = 0\}} \right]   = \\  \left(\frac{2-\g}{\g-1}\right)^{q/2} ((\gamma-1)Ch)^{\frac{\gamma-2}{2(\gamma-1)}q}
\E\left[\left[\left(\frac{V_0^{1-\gamma}}{(\gamma-1)Ch}\right)^\frac{\gamma-2}{\gamma-1}- \left(\frac{V_0^{1-\gamma}}{(\gamma-1)Ch}+1\right)^\frac{\gamma-2}{\gamma-1}\right]^{q/2}\right] e^{-\l h}. 
\end{multline}
To estimate this last expression we need, letting $a := \frac{2-\g}{\g-1}$, the following modifications of \eqref{tecineq1}, valid for $z>0$:
\be{tecineq2}
\begin{array}{ll}
a(z+1)^{-1} z^{-a} \leq z^{-a} - (z+1)^{-a}  \leq (z+1)^{-1} z^{-a} & \mbox{ for } 0<a \leq 1 \\
(z+1)^{-1} z^{-a} \leq z^{-a} - (z+1)^{-a}  \leq a(z+1)^{-1} z^{-a} & \mbox{ for } a>1.
\end{array}
\ee
Using these inequalities as in \eqref{main2} we obtain, for some $C>1$
\be{main21}
\frac{1}{C} \E\left[  \left(\frac{V_0^{\gamma-1} h }{1+V_0^{\gamma-1} h} V_0^{2-\gamma}\right)^q \right] \leq \E\left[\left(\int_0^h V_t dt\right)^{q/2} \ind_{\{i(h) = 0\}} \right] \leq C \E\left[  \left(\frac{V_0^{\gamma-1} h }{1+V_0^{\gamma-1} h} V_0^{2-\gamma}\right)^q \right].
\ee
We now observe that
\be{main22}
\begin{split}
\E\left[  \left(\frac{V_0^{\gamma-1} h }{1+V_0^{\gamma-1} h} V_0^{2-\gamma}\right)^{q/2} \right]   =   & \ \E\left[  \left(\frac{V_0^{\gamma-1} h }{1+V_0^{\gamma-1} h} V_0^{2-\gamma}\right)^{q/2} \ind_{\{V_0^{\gamma-1} h \leq 1\}} \right] \\ & + \E\left[  \left(\frac{V_0^{\gamma-1} h }{1+V_0^{\gamma-1} h} V_0^{2-\gamma}\right)^{q/2} \ind_{\{V_0^{\gamma-1} h > 1\}} \right] \\ \sim & \ h^{q/2}  \E\left[ V_0^{q/2} \ind_{\{V_0^{\gamma-1} h \leq 1\}} \right] + \E\left[ V_0^{\frac{q}{2}(2-\g)} \ind_{\{V_0^{\gamma-1} h > 1\}} \right] .
\end{split}
\ee
In order to estimate the two summand of the left hand side of \eqref{main22} we use the following fact, whose simple proof follows from simple point wise bounds, and it is omitted.
Let $\mu$ be a probability on $[0,+\infty)$ such that $\mu((u,+\infty))$ is regularly varying with exponent $-\xi<0$. Then 
\begin{eqnarray} 
\int_0^x u^p \mu(du) \sim x^{p-\xi} & & \mbox{for } p>\xi \label{main23} \\
\int_x^{+\infty}  u^p \mu(du) \sim x^{p-\xi} & & \mbox{for } p<\xi . \label{main24}
\end{eqnarray}
Let $\mu$ be the law of $V_0$, so that, by Theorem \ref{th:stmeas}, $\xi = \a + \g -1$. Since $q>2( \a + \g -1)$, by \eqref{main23} we have $\E\left[ V_0^{q/2} \ind_{\{V_0^{\gamma-1} h \leq 1\}} \right] \sim h^{-\frac{1}{\g-1}(\frac{q}{2}-\a-\g+1)}$, and therefore
\be{main25}
h^{q/2}  \E\left[ V_0^{q/2} \ind_{\{V_0^{\gamma-1} h \leq 1\}} \right] \sim h^{\frac{\gamma-2}{2(\gamma-1)}q+ \frac{\a+\g-1}{\g-1}}.
\ee
Moreover, by \eqref{main24}, also
\be{main26}
\E\left[ V_0^{\frac{q}{2}(2-\g)} \ind_{\{V_0^{\gamma-1} h > 1\}} \right] \sim h^{\frac{\gamma-2}{2(\gamma-1)}q+ \frac{\a+\g-1}{\g-1}}.
\ee
for $\frac{q}{2}(2-\g) < \a+\g-1$, while
\be{main27}
\E\left[ V_0^{\frac{q}{2}(2-\g)} \ind_{\{V_0^{\gamma-1} h > 1\}} \right] = +\infty
\ee
for $\frac{q}{2}(2-\g) > \a+\g-1$.

Summing up, we have shown that
\begin{eqnarray}
\E\left[\left(\int_0^h V_t dt\right)^{q/2} \ind_{\{i(h) = 0\}} \right]  \sim & h^{\frac{\gamma-2}{2(\gamma-1)}q+ \frac{\a+\g-1}{\g-1}} &  \mbox{ for } 2(\a+\g-1) < q < \frac{2(\a+\g-1)}{2-\g}  \nonumber \\ \E\left[\left(\int_0^h V_t dt\right)^{q/2} \ind_{\{i(h) = 0\}} \right] 
= & +\infty &   \mbox{ for } q > \frac{2(\a+\g-1)}{2-\g}. \label{main271}
\end{eqnarray}

\noindent
{\em Case $i(h) = 1$}. This case is dealt with as for $\g>2$, and, using the same argument leading to \eqref{main21}, one sees that the crucial term to estimate is
\be{main28}
\E\left[  \left(\frac{V_{T_1}^{\gamma-1} (h-T_1) }{1+V_{T_1}^{\gamma-1} (h-T_1)} V_{T_1}^{2-\gamma}\right)^{q/2}  \ind_{\{i(h) = 1\}}\right].
\ee
Since $V_{T_1} \geq L_{T_1}$, \eqref{main28} can be bounded from below by
\[
\E\left[  \left(\frac{L_{T_1}^{\gamma-1} (h-T_1) }{1+L_{T_1}^{\gamma-1} (h-T_1)} L_{T_1}^{2-\gamma}\right)^{q/2}  \ind_{\{i(h) = 1\}}\right]
\]
which takes the value infinity as soon as $\E\left[ \left(L_{T_1}^{2-\gamma}\right)^{q/2}  \ind_{\{i(h) = 1\}}\right] = +\infty$. Recalling that $L_{T_1}$ independent of $\{i(h) =1 \}$ and it has law $\nu$, this hold as $q>\frac{2\a}{2-\g}$. This implies that
\be{main281}
\E\left[\left(\int_0^h V_t dt\right)^{q/2} \ind_{\{i(h) = 1\}} \right] = +\infty \ \ \ \mbox{ for } q>\frac{2\a}{2-\g}.
\ee
Comparing with \eqref{main271}, note that $\frac{2\a}{2-\g} < \frac{2(\a+\g-1)}{2-\g}$. Thus, assume $2(\a+\g-1)< q <\frac{2\a}{2-\g}$ (note that, being by assumption $\a + \g>2$, indeed $2(\a+\g-1) < \frac{2\a}{2-\g}$). An upper bound for \eqref{main28} is given by
\begin{multline} \label{main29}
 \E\left[  \left(\frac{V_{T_1}^{\gamma-1} (h-T_1) }{1+V_{T_1}^{\gamma-1} (h-T_1)} V_{T_1}^{2-\gamma}\right)^{q/2}  \ind_{\{i(h) = 1\}}\right]  \\  \leq \E\left[  \left(\frac{(V_0 + L_{T_1})^{\gamma-1} h }{1+(V_0 + L_{T_1})^{\gamma-1} h} (V_0 + L_{T_1})^{2-\gamma}\right)^{q/2}  \ind_{\{i(h) = 1\}}\right] \\=  \E\left[  \left(\frac{(V_0 + L_{T_1})^{\gamma-1} h }{1+(V_0 + L_{T_1})^{\gamma-1} h} (V_0 + L_{T_1})^{2-\gamma}\right)^{q/2} \right] \P(i(h) = 1),
\end{multline}
where we used the facts that $V_0$ and $L_{T_1}$ are independent of $\{i(h) = 1\}$. Now, \eqref{main29} is estimated exactly as \eqref{main22}, but with $V_0 + L_{T_1}$ in place of $V_0$. Since the tails of $V_0 + L_{T_1}$ are the same of those of $L_{T_1}$, i.e. regularly varying with exponent $\a$, while $\P(i(h)=1) \sim h$, we get
\[
\E\left[  \left(\frac{(V_0 + L_{T_1})^{\gamma-1} h }{1+(V_0 + L_{T_1})^{\gamma-1} h} (V_0 + L_{T_1})^{2-\gamma}\right)^{q/2} \right]  \sim h^{\frac{\gamma-2}{2(\gamma-1)}q+ \frac{\a}{\g-1}} \P(i(h) = 1) \sim h^{\frac{\gamma-2}{2(\gamma-1)}q+ \frac{\a + \g -1}{\g-1}}.
\]
Summing up:
\begin{eqnarray}
\E\left[\left(\int_0^h V_t dt\right)^{q/2} \ind_{\{i(h) = 1\}} \right]  \sim & h^{\frac{\gamma-2}{2(\gamma-1)}q+ \frac{\a+\g-1}{\g-1}} &  \mbox{ for } 2(\a+\g-1) < q < \frac{2\a}{2-\g}  \nonumber \\ \E\left[\left(\int_0^h V_t dt\right)^{q/2} \ind_{\{i(h) = 0\}} \right] 
= & +\infty &   \mbox{ for } q > \frac{2\a}{2-\g} . \label{main30}
\end{eqnarray}

\noindent
{\em Case $i(h) \geq 2$}. This case goes along the same line as for $\g>2$, using the upper bound obtained for $i(h) = 1$. The details are omitted. The proof for $1<\g<2$ is thus completed.

\noindent
\textbf{Part 3: $\g=2$}\\

In this case we have, in the case of no jumps ($i(h) = 0$),
\begin{equation*}
V_t=\left( V_0^{-1}+Ct\right)^{-1}
\end{equation*}
and therefore
\be{main31}
\int_0^h V_t dt=\frac{1}{C}\left[\log(V_0^{-1}+Ch)-\log(V_0^{-1})\right]=
\frac{1}{C}\left[\log(1+ChV_0)\right].
\ee
Un upper bound for $\E\left[\left(\int_0^h V_t dt\right)^{q/2} \ind_{\{i(h) = 0\}} \right]$ is obtained using \eqref{main31} and  the inequality, valid for $y,r>0$,
\[
\log(1+y) \leq \frac{1}{r} y^r,
\]
which gives
\[
\E\left[\left(\int_0^h V_t dt\right)^{q/2} \ind_{\{i(h) = 0\}} \right] \leq \frac{1}{C^{q/2}} C^{rq/2} h^{rq/2} \E\left(V_0^{rq/2} \right).
\]
Since $\E\left(V_0^{rq/2} \right) < +\infty$ for $\frac{rq}{2} < \a+1$, letting $\frac{rq}{2} \uparrow \a+1$ we obtain
\be{main32}
\E\left[\left(\int_0^h V_t dt\right)^{q/2} \ind_{\{i(h) = 0\}} \right] \leq C h^{rq/2},
\ee
for some $C>0$ and every $r$ such that $\frac{rq}{2} < \a+1$.
A corresponding lower bound is obtained using the inequality 
\[
\log(1+y) \geq \frac{1}{2} \ind_{(1,+\infty)}(y),
\]
which gives
\be{main33}
\E\left[\left(\int_0^h V_t dt\right)^{q/2} \ind_{\{i(h) = 0\}} \right] \geq \frac{1}{(2C)^{q/2}} \P(Ch V_0 > 1) \sim h^{\a+1},
\ee
where we have used Theorem \ref{th:stmeas} for the last inequality. By \eqref{main32} and $\eqref{main33}$ we have
\[
\E\left[\left(\int_0^h V_t dt\right)^{q/2} \ind_{\{i(h) = 0\}} \right]  \sim h^{\a+1}.
\]

The cases with $i(h) \geq 1$ are similar to what seen in Parts 1 and 2, and are omitted.

\epr

\bigskip

\noindent
{\em Proof of Theorem \ref{th:main}}. We now complete the proof of Theorem \ref{th:main}. We need to extend Proposition \ref{prop:main1} in two directions: a) generalize from $f(v) = Cv^{\g}$ to any $f$ satisfying Assumption {\bf A4}; b) extend to Levy subordinator satisfying Assumptions {\bf A2} and {\bf A3}, thus with a possibly infinite characteristic measure $\nu$.

\noindent
{\em Step 1}. We keep all assumption of Proposition \ref{prop:main1}, except that we require $f(v) = Cv^{\g}$ only for $v>\e$, for some $\e>0$, and $f$ satisfies Assumption {\bf A4}. In other words we do not prescribe the asymptotics of $f$ near $v=0$. Let $V,V'$ be solutions, respectively, of the equations
\[
\begin{split}
dV_t & = \, - f(V_t) dt + dL_t \\ 
dV'_t & = \, -CV_t'^{\g} + dL_t.
\end{split}
\]
Assume $V_0 = V'_0 = v > 0$. We claim that 
\be{th-1}
|V_t - V'_t| \leq 2\e 
\ee
 a.s., for every $t \geq 0$. This follows from the following fact: there is a constant $\d>0$ such that as soon as  $|V_t - V'_t| \geq 2\e$,
\be{th-2}
\frac{d}{dt}  |V_t - V'_t|  \leq -\d.
\ee
To see \eqref{th-2}, suppose first $V_t - V'_t \geq 2\e$. In particular $V_t \geq \e$, so
\[
\frac{d}{dt} [V_t - V'_t] = -C\left(V_t^{\g} - V_t'^{\g} \right) < -C(2\e)^{\g},
\]
where we have used the fact that, for $c>0$, the map $(x+c)^{\g} - x^{\g}$ is increasing for $x>0$. Suppose now $V'_t - V_t \geq 2\e$. If $V_t \geq \e$ then,
\[
\frac{d}{dt} [V'_t - V_t] = -C\left(V_t'^{\g} - V_t^{\g} \right) < -C(2\e)^{\g};
\]
If $V_t < \e$, since $f$ is increasing, then
\[
\frac{d}{dt} [V'_t - V_t] = -C V_t'^{\g} + f(V_t) \leq -C(2\e)^{\g} +C\e^{\g} < 0.
\]
Thus \eqref{th-2}, and so \eqref{th-1} is proved. In particular, the law of $V_t$ is stochastically smaller than that of $V'_t + 2\e$, which means that for every $g$ increasing and bounded, $E[g(V_t)] \leq E[g(V'_t + 2\e)]$. By the ergodicity results proved in \cite{sam}, this inequality can be taken to the limit as $t \ra +\infty$, so to a stochastic inequality between the stationary distributions of $V$ and $V'$. This implies that we can realize, on a suitable probability space, two random variables $V_0$ and $V'_0$, independent of the Levy process $L$, distributed according to the stationary laws of the corresponding processes, and such that $V_0 \leq V'_0 + 2\e$. By repeating the argument above, we see that the inequality $V_t \leq V'_t + 2\e$ is a.s. preserved for all $t>0$ also for the stationary processes. It follows that
\be{th-3}
\E\left[\left(\int_0^h V_t dt\right)^{q/2}\right] \leq \E\left[\left(\int_0^h [V'_t+2\e] dt\right)^{q/2}\right] \leq 2^{q/2-1} \left\{\E\left[\left(\int_0^h V'_t dt\right)^{q/2}\right] + (2\e h)^{q/2} \right\}.
\ee
Since
\be{ext0}
A(q) = \limsup_{h \ra 0} \frac{\log \E\left[\left(\int_0^h V_t dt\right)^{q/2}\right]}{\log h},
\ee
and
\be{ext1}
A'(q) = \lim_{h \ra 0} \frac{\log \E\left[\left(\int_0^h V'_t dt\right)^{q/2}\right]}{\log h} \leq \frac{q}{2},
\ee
by \eqref{th-3} we get
\[
A(q) \geq A'(q).
\]
By exchanging the role of $V$ and $V'$ we get $A(q) = A'(q)$. Moreover, the existence of the limit \eqref{ext1}, which follows from Proposition \ref{prop:main1}, implies that also \eqref{ext0} is a limit.  Since $A'(q)$ is given by Proposition \ref{prop:main1}, the first extension is obtained.

\noindent
{\em Step 2}. In this step we allow the Levy process $L$ to have infinite characteristic measure $\nu$ and {\em positive drift} $m>0$, though satisfying Assumptions {\bf A2} and {\bf A3}. On the other hand we make a specific choice for $f$:  $f(v) = Cv^{\g}$ for $v>\e$, while $f$ is linear in $[0,\e)$, with $f(0) = 0$ and $f(\e) = C \e^{\g}$. Moreover we let $\nu_{\e} := \nu \ind_{[\e,+\infty)}$, which is a finite measure. Denote by $L^{(\e)}$ the compound Poisson process with characteristic measure $\nu^{\e}$,
and by $V^{(\e)}$ the solution of 
\be{th-4}
dV^{(\e)}_t  = \, - f(V^{(\e)}_t) dt + dL^{(\e)}_t
\ee
The original Levy process $L$ can be decomposed in the form $L_t = L^{(\e)} + L^{(<\e)}$, where $L^{(<\e)}$ is independent of $L^{(\e)}$, it has characteristic measure $\nu_{\e} := \nu \ind_{[0,\e)}$ and drift $m>0$. Writing
\be{th-5}
dV_t  = \, - f(V_t) dt + dL_t,
\ee
we obtain
\be{th-6}
d(V_t - V^{(\e)}_t) = \, - [f(V_t)- f(V^{(\e)}_t)]dt + dL^{(<\e)}.
\ee
This implies, for instance that whenever $V^{(\e)}_0 \leq V_0$, then $V^{(\e)}_t \leq V_t$ for all $t>0$. Thus, using as above the ergodicity of $V$ and $V^{(\e)}$, $ V_t$ dominates stochastically $V^{(\e)}_t$ also in equilibrium. Thus, as before, we can start the processes from $V^{(\e)}_0$ and $V_0$, each having the corresponding stationary distribution, and such that $V^{(\e)}_0 \leq V_0$. Thus $V^{(\e)}_t \leq V_t$ for all $t>0$. Note that with this construction we have that the two processes in \eqref{th-4} and \eqref{th-5} are separately stationary, by not necessarily the Markov process $(V^{(\e)}_t , V_t)$, whose law will be demoted by $\mu_t^{(2)}$, is stationary. To fix this we observe that, since the family of distribution $(\mu_t^{(2)})_{t \geq 0}$ is tight, by a standard argument its Cesaro means $\frac{1}{t} \int_0^t \mu_s^{(2)} ds$ admit at least a limit point, which is a stationary distribution for $(V^{(\e)}_t , V_t)$. This limiting operation preserves the stochastic order between the laws of the two components. Thus, we can assume to realize $V^{(\e)}_0$ and $V_0$ in such a way their {\em joint} distribution is stationary for \eqref{th-4} and \eqref{th-5}, and $V^{(\e)}_0 \leq V_0$.

Now we use the fact that $f$ is super linearly increasing, to conclude that
\[
f(V_t)- f(V^{(\e)}_t) \geq c[V_t - V^{(\e)}_t]
\]
for some $c>0$. It follows that
\[
d(V_t - V^{(\e)}_t)  \leq \, - c[V_t - V^{(\e)}_t] + dL^{(<\e)}_t, 
\]
which implies that
\be{th-7}
0 \leq V_t - V^{(\e)}_t \leq e^{-ct} [V_0 - V^{(\e)}_0] + \int_0^t e^{-c(t-s)} dL^{(<\e)}_s.
\ee
Since the law of $V_t - V^{(\e)}_t$ does not depend on $t$, it must be stochastically dominated by the limit of the law of the r.h.s. of \eqref{th-7}, which is just the stationary distribution of the Ornstein-Uhlenbeck process
\[
dZ_t = -cZ_t dt + dL^{(<\e)}_t.
\]
As observed e.g. in \cite{kl1}, this stationary law is infinitely divisible with characteristic pair $(m,\tilde{\nu})$, with
\[
\tilde{\nu}([x,+\infty)) = \int_{[x,+\infty)} u^{-1} \nu_{\e}(du).
\]
Since $\nu_{\e}$, and therefore $\tilde{\nu}$, has bounded support, the stationary law of $Z_t$ has moments of all order (see e.g. \cite{wol}).
So, also $V_t - V^{(\e)}_t$ has moments of all order. Thus, using the inequality $(x+y)^q \leq 2^{q-1}[x^q + y^q]$ for $x,y \geq 0$, we have
\begin{multline*}
\E\left[\left( \int_0^h V^{(\e)}_t dt \right)^{q/2} \right] \leq \E\left[\left( \int_0^h V_t dt \right)^{q/2} \right]  \\ \leq 2^{q/2-1} \left\{\E\left[\left( \int_0^h V^{(\e)}_t dt \right)^{q/2} \right]  + \E\left[\left( \int_0^h [V_t - V^{(\e)}_t]dt \right)^{q/2} \right] \right\} \\  \leq 2^{q/2-1} \left\{\E\left[\left( \int_0^h V^{(\e)}_t dt \right)^{q/2} \right]  + h^{q/2}\E\left[\left(V_0 - V^{(\e)}_0 \right)^{q/2} \right] \right\}. 
\end{multline*}
Since, by Proposition \ref{prop:main1} and steps 1, $\E\left[\left( \int_0^h V^{(\e)}_t dt \right)^{q/2} \right] \sim h^{A(q)}$ for every $\e>0$ and $A(q) \leq q/2$, it follows that
\[
\E\left[\left( \int_0^h V_t dt \right)^{q/2} \right]  \sim h^{A(q)},
\]
thus completing the proof of this step.

\noindent
{\em Step 3}. The extension of Proposition \ref{prop:main1} to any $f$ which satisfy ({\bf A4}) is now easy, and it will only be sketched. In a first stage, repeating the argument in step 1, one extends from the special $f$'s used for step 2, to the larger class of $f$ in step 1. 

The further extension to a general $f$ which satisfy ({\bf A4}) proceeds as follows:
for every $\d>0$ we can find $f_1$ and $f_2$ such that $f_1 \leq f \leq f_2$, and $f_1(v) = C_1 v^{\g - \d}$, $f_2(v) = C_2 v^{\g +\d}$ for $v>\e$. By using coupling arguments similar to those in step 1, one shows that the scaling function $A(q)$ of the process with drift $f$ is bounded above and below by the scaling functions of the processes with drift $f_1$ and $f_2$. The continuity of $A(q)$ w.r.t. $\g$, and the fact that $\d$ is arbitrary, implies that $A(q)$ is given by Proposition \ref{prop:main1}.

\qed

\noindent \textbf{\textit {Acknowledgments:}} We thank Francesco Caravenna for useful discussions and suggestions.

\end{document}